 \documentclass[journal,final]{IEEEtran}

\usepackage[pdftex]{graphicx}
\usepackage{amsfonts}
\usepackage{setspace}
\usepackage{wrapfig}
\usepackage{color}
\usepackage{cancel}
\usepackage{url}
\usepackage{enumerate}
\usepackage{cite}
\usepackage{amssymb}
\usepackage[cmex10]{amsmath}
\usepackage{subfigure}
\def\figurename{Fig.}

\newcommand{\field}[1]{\mathbb{#1}}                                   
\newcommand{\prn}[1]{\left(#1\right)}                                 
\newcommand{\bigprn}[1]{\Bigl(#1\Bigr)}                
\newcommand{\Prn}[1]{\left[#1\right]}                                 
\newcommand{\dint}[3]{\displaystyle\int_{#1}^{#2}{#3}}                
\newcommand{\abs}[1]{\left|#1\right|}                                 
\newcommand{\eltwo}[1]{\left\|#1\right\|_2}                           
\newcommand{\elinf}[1]{\left\|#1\right\|_\infty}                      
\newcommand{\norm}[1]{\left\|#1\right\|}                              
\newcommand{\sign}[1]{\ \mathrm{sign}{\left(#1\right)}}								
\newcommand{\set}[1]{\displaystyle\left\{#1\right\}}									
\newcommand{\E}[1]{\mathbb{E}\set{#1}}									

\DeclareMathOperator*{\argmin}{arg\,min}															
\newtheorem{theorem}{Theorem}
\newtheorem{lemma}{Lemma}

\newtheorem{definition}{Definition}

\newcommand{\G}{\Gamma}
\newcommand{\Gk}{\Gamma_k}
\newcommand{\Gc}{\Gamma^c}

\newcommand{\PP}{\Phi_{\Gamma}^T\Phi_{\Gamma}}
\newcommand{\PPk}{\Phi_{\Gamma_k}^T\Phi_{\Gamma_k}}

\newcommand{\PcP}{\Phi_{\Gamma^c}^T\Phi_{\Gamma}}

\newcommand{\GD}{\Delta}

\newcommand{\adag}{a^{\dagger}}

\usepackage{scrtime} 
\usepackage[draft,scrtime]{prelim2e}

\begin{document}
%
\title{Convergence Speed of a Dynamical System for Sparse Recovery}

%
%
%

\author{Aur\`ele~Balavoine,~\IEEEmembership{Student Member,~IEEE,}
		Christopher~J.~Rozell,~\IEEEmembership{Senior Member,~IEEE,}
        and~Justin~Romberg,~\IEEEmembership{Senior Member,~IEEE,}
\thanks{The authors are with the School
of Electrical and Computer Engineering, Georgia Institute of Technology, Atlanta,
GA, 30332-0250 USA}\thanks{Email: \{aurele.balavoine,crozell,jrom\}@gatech.edu.}}

\maketitle

\begin{abstract}

This paper studies the convergence rate of a continuous-time dynamical system for $\ell_1$-minimization, known as the Locally Competitive Algorithm (LCA). Solving \mbox{$\ell_1$-minimization} problems efficiently and rapidly is of great interest to the signal processing community, as these programs have been shown to recover sparse solutions to underdetermined systems of linear equations and come with strong performance guarantees.  The LCA under study differs from the typical $\ell_1$-solver in that it operates in continuous time: instead of being specified by discrete iterations, it evolves according to a system of nonlinear ordinary differential equations.  The LCA is constructed from simple components, giving it the potential to be implemented as a large-scale analog circuit.

The goal of this paper is to give guarantees on the convergence time of the LCA system. To do so, we analyze how the LCA evolves as it is recovering a sparse signal from underdetermined measurements.  We show that under appropriate conditions on the measurement matrix and the problem parameters, the path the LCA follows can be described as a sequence of linear differential equations, each with a small number of active variables.  This allows us to relate the convergence time of the system to the restricted isometry constant of the matrix. Interesting parallels to sparse-recovery digital solvers emerge from this study. Our analysis covers both the noisy and noiseless settings and is supported by simulation results.
\end{abstract}

\makeatletter{\renewcommand*{\@makefnmark}{}
\footnotetext{This work was partially supported by NSF grant CCF-0905346.}\makeatother}

\makeatletter{\renewcommand*{\@makefnmark}{}
\footnotetext{Submitted to the IEEE Transactions on Signal Processing on December 13, 2012.}\makeatother}

\section{Introduction}
\label{sec:intro}

\PARstart{C}{ompressed} Sensing (CS) has triggered extensive research because of compelling results on the reconstruction of sparse signals (i.e. signals with few non-zero elements) from highly-undersampled linear measurements.  The main results of CS show that coded measurements can be used to simultaneously acquire and compress a signal, requiring many fewer resources (e.g., time, storage, etc.) than traditional sampling approaches.
However, the process of reconstructing the original signal from its compressed measurements requires a significant amount of computation and remains a bottleneck in the processing pipeline. 

The approach to signal reconstruction that has been most extensively studied involves 
solving an optimization program that minimizes a combination of a mean-squared error term and a sparsity-inducing term (typically measured using the $\ell_1$-norm).  
Specifically, given a set of (possibly noisy) measurements $y\in\mathbb{R}^M$ of a signal $\adag\in\mathbb{R}^N$ through a $M\times N$ matrix $\Phi$, we estimate $\adag$ by solving
\begin{equation}
\hat{a}^{\dagger} = \argmin\limits_{a}{\dfrac{1}{2}\left\|y-\Phi a\right\|^2_2 + \lambda \left\|a\right\|_1}, 
\label{eq:l1}
\end{equation}
where $\norm{a}_1 = \sum_i \abs{a_i}$.
Despite the optimization in~\eqref{eq:l1} being a convex and tractable program with many specialized solvers (e.g., \cite{donoho_fast_2008,chen_atomic_2001,kim_interior-point_2007,figueiredo_gradient_2007,becker_nesta:_2011,yin_bregman_2008,daubechies_iterative_2004}), the required computation in most problem sizes of interest makes it  prohibitive to perform CS reconstruction in real time or on low-power embedded platforms.

The Locally Competitive Algorithm (LCA)~\cite{rozell_sparse_2008} is a \emph{continuous-time} system of coupled nonlinear differential equations that settles to the minimizer of \eqref{eq:l1} in steady state~\cite{balavoine_convergence_2012}.  The LCA architecture consists of simple components (matrix-vector operations and a pointwise nonlinearity for thresholding), giving it the potential to be implemented in an analog circuit~\cite{shapero_low_2012,shapero_configurable_2013}.  Analog networks for solving optimization problems have a long history, dating back to Hopfield's pioneering results for linear programming~\cite{hopfield_neural_1982} (a comprehensive treatment of the subject can be found in~\cite{cichocki_neural_1993}).  Such analog systems can potentially have significant speed and power advantages over their digital counterparts.

While analog implementations of the LCA have the potential to alleviate the bottleneck of CS signal reconstruction in some scenarios, as with any signal processing system, it is important to have strong performance guarantees before deploying the system in an application.
Prior work~\cite{balavoine_convergence_2012} has studied the convergence behavior of the LCA in a general setting (with no assumption on the signal or the matrix) to prove that the system has global asymptotic convergence to the correct solution.  In this general setting, the LCA is shown to converge exponentially fast, provided some condition that depends on the solution path of the system (i.e., which nodes cross threshold to become active during convergence). More specifically, as the nodes evolve, the LCA dynamics switch between sets of linear ordinary differential equations that involve submatrices of $\Phi$. If each submatrix is well-conditioned, then the exponential convergence result follows. The main contribution of this paper is to study the specific case of sparse recovery. Interestingly, our analysis depends on the well-known Restricted Isometry Property (RIP) for CS measurement matrices. This condition ensures that every submatrix of a specific size is well-conditioned. The results in this paper establish conditions on problem parameters (such as signal sparsity $S$, ambient dimension $N$, and number of measurements $M$) that guarantee that the size of each submatrix is indeed small. These guarantees can then be used to provide strong bounds on the convergence speed of the system. Our resulting conditions are naturally analogous to existing bounds for traditional digital algorithms.

After reviewing the guarantees for existing digital algorithms and prior analysis of the LCA in Section~\ref{sec:algo}, we present our two main results in Section~\ref{sec:charac}.  Theorem~\ref{th:opt} establishes conditions ensuring that only nodes that are part of the support of the original signal $\adag$ become active during convergence.
Theorem~\ref{th:size} relaxes these conditions to allow a fixed number of nodes to enter the support during convergence. 
Both results allow us to establish a bound on the exponential rate of convergence for the LCA by relating it to the restricted isometry constant of the matrix $\Phi$. Section~\ref{ssec:discussion} explores the implication of the two theorems when the measurement matrix is random. The qualitative predictions of these theoretical guarantees are explored in simulation in Section~\ref{sec:sim}.

\section{Background and Related work}
\label{sec:algo}

The analysis in this paper differs from previous studies that appear in the CS literature because of the continuous nature of the LCA algorithm. In particular, the complexity of the LCA cannot be expressed in terms of a number of ``iterations'' as is often done for digital algorithms. Nevertheless, some analogies to previous work can be drawn. In this section, we describe some existing approaches to sparse recovery, their associated guarantees, and how they relate to the LCA.

\subsection{Existing algorithms}
\label{ssec:algo}
\subsubsection{$\ell_1$-solvers}

The $\ell_1$-norm in \eqref{eq:l1} is used as a surrogate for the ideal pseudo-norm $\norm{a}_0$, which counts the number of non-zero elements.  Under particular conditions on $\Phi$, it can be shown that the performance of the relaxed program \eqref{eq:l1} is comparable to the idealized (but generally intractable~\cite{natarajan_sparse_1995}) sparse approximation problem.
One such condition is known as the Restricted Isometry Property (RIP). The RIP guarantees that every submatrix formed from a small subset of columns of $\Phi$ is a near isometry.
\begin{definition}
\label{def:rip}
The matrix $\Phi$ satisfies the RIP of order $k$ if there exist a constant $\delta\in\prn{0,1}$, such that for any vector $x\in\field{R}^N$ such that $\norm{x}_0\leq k$, we have:
\begin{equation}
(1-\delta)\eltwo{x}^2\leq\eltwo{\Phi x}^2\leq(1+\delta)\eltwo{x}^2.
\label{eq:rip}
\end{equation}
We also say that $\Phi$ satisfies the RIP with parameters $(k,\delta)$. The RIP-constant $\delta_k$ of order $k$ for $\Phi$ is defined as the smallest positive constant $\delta$ satisfying $\eqref{eq:rip}$.
\end{definition}
The main advantage of solving \eqref{eq:l1} is the existence of sharp results on the $\ell_2$- and $\ell_1$-norm of the error~\cite{candes_restricted_2008}.
%
However, despite the many digital algorithms proposed, the complexity involved with solving \eqref{eq:l1} remains a challenge. State-of-the art solvers can handle large scale problems, but lack strong guarantees about their running time~\cite{kim_interior-point_2007,figueiredo_gradient_2007,van_den_berg_probing_2008,becker_nesta:_2011,yin_bregman_2008}. On the other hand, iterative thresholding schemes are simple and come with guarantees about the number of iterations needed to achieve a certain accuracy~\cite{daubechies_iterative_2004,blumensath_iterative_2008}, but they may require many computationally expensive iterations. Homotopy-based schemes find the solution to \eqref{eq:l1} by tracing a piecewise-linear solution path as the tradeoff parameter $\lambda$ is varied.  If the solution is very sparse, these approaches can converge in exactly $S$ iterations (known as the \textit{S-step property})~\cite{donoho_fast_2008}. 

The LCA differential equation resembles a continuous-time version of an iterative thresholding step. However, because the system evolves according to a piecewise-linear dynamical system that changes each time a node crosses threshold, its solution path is very similar to the Homotopy~\cite{donoho_fast_2008} and its approximate version LARS~\cite{efron_least_2004}. Much in the spirit of the S-step property and for a similar number of measurements for random matrices, Theorem~\ref{th:opt} shows that the LCA active set remains a subset of the $S$ optimal nodes.

\subsubsection{Greedy Algorithms}
A second approach to sparse signal recovery is through the use of iterative greedy algorithms.  A common solver in this family is Orthogonal Matching Pursuit (OMP), which at each iteration adds to the support the element that has the strongest correlation with the residual. Through the use of the RIP, OMP was shown to recover an $S$-sparse signal in exactly $S$ iterations (i.e., has the $S$-step property) in the noiseless setting~\cite{davenport_analysis_2010}, for a number of random measurements similar to our first result in Theorem~\ref{th:opt}.

Letting OMP run for more than $S$ iterations, recovery results have been obtained in the noiseless~\cite{tropp_signal_2007} and noisy cases~\cite{zhang_sparse_2011} for fewer random measurements. Similar recovery results exist for Regularized Orthogonal Matching Pursuit (ROMP)~\cite{needell_uniform_2009} and Compressive Sampling Matching Pursuit (CoSaMP)~\cite{needell_cosamp:_2008}. In contrast to OMP, ROMP and CoSaMP add a set of nodes at each iteration. While in the case of the LCA, we are not concerned with iteration count, our Theorem~\ref{th:size} is similar in nature. Letting the active set grow larger than the $S$ optimal nodes still yield guarantees on the convergence time, while reducing the number of CS measurements necessary.


\subsection{The Locally Competitive Algorithm}
\label{sec:lca}

\subsubsection{LCA structure and dynamics}
\label{ssec:lca}

The LCA can be viewed as a network of nodes that evolve according to a set of coupled, nonlinear Ordinary Differential Equations (ODEs):
\begin{equation} 
\begin{split}
	\tau \dot{u}(t)&=-u(t) - (\Phi^T \Phi - I) \ a(t) + \Phi^Ty\\
	a(t) &= T_{\lambda}(u(t))\\
\end{split} .
\label{eq:dyn}
\end{equation}
These dynamics govern the set of internal state variables, $u_n(t)$ for $n=1, \ldots, N$,  each associated with a single dictionary element $\Phi_n$. The internal states produce output variables $a_n(t)$ for $n=1, \ldots, N$ through a nonlinear pointwise activation function $T_{\lambda}(\cdot)$. 
The time constant $\tau$ is determined by the physical properties of the solver implementing the system. Since $\tau$ does not affect the mathematical analyses of the system, we often take $\tau=1$ except when we want to stress its influence on the convergence speed. We assume throughout that the columns of $\Phi =\left[\Phi_1, \ldots, \Phi_N\right]$ have unit norm: $\norm{\Phi_n}_2=1$. 
%
%
To solve \eqref{eq:l1}, the activation function used is the soft-thresholding function~\cite{rozell_sparse_2008}:
\begin{equation}
a_n(t) = T_{\lambda}(u_n(t)) = \begin{cases}
0, & \abs{u_n(t)} \leq \lambda \\
u_n(t)-\lambda z_n(t), & \abs{u_n(t)} > \lambda 
\end{cases},
\label{eq:thresh}
\end{equation}
where $z_n(t) = \sign{u_n(t)}$ is the sign of the $n^{th}$ internal state variable.
Though $\ell_1$-minimization is our focus, recent work has shown that many other sparsity-inducing penalty functions can be minimized in the same system by changing the form of $T_{\lambda}(\cdot)$~\cite{charles_common_2012}.

\subsubsection{Notations}
It can be seen from \eqref{eq:thresh} that the activation function is composed of two operating regions.  When \mbox{$\abs{u_n}\leq\lambda$}, the output $a_n$ is zero and we call the node \emph{inactive}. When  $\abs{u_n} > \lambda$, the output $a_n$ is strictly increasing with $u_n$ and we call the node \emph{active}.  Denote by $\G(t)$ the current \emph{active set} (i.e., the set of indices $\G(t)=\{k\in\Prn{1,N},\ \abs{u_k(t)}>\lambda\}$), and denote by $\Gc(t)$ the \emph{inactive set} consisting of nodes that are below threshold. While the active set changes with time as the network evolves, for the sake of readability and when it is clear from the context, we omit the dependence on time in the notation and just write the active set as $\G$. The sequence of \emph{switching times} for which the system moves from the set of active nodes $\G_{k-1}$ to $\G_{k}$ is the sequence $\{t_k\}_{\{k\in\field{N}\}}$. In the following, we also denote by $\Phi_{\mathcal{T}}$ the matrix composed of the columns of $\Phi$ indexed by the set $\mathcal{T}$, setting all the other entries to zero. Similarly, $u_\mathcal{T}$ and $a_\mathcal{T}$ refer to the elements in the original vectors indexed by $\mathcal{T}$ setting other entries to zero. 

\subsubsection{LCA Convergence Speed}

As with their digital counterparts, it is desirable to know how fast continuous-time systems such as the LCA converge. 
%
The LCA has been shown to be exponentially convergent~\cite{balavoine_convergence_2012}, with a convergence speed that depends on the transient activity in the system.  To state this result,  define $\G_*$ as the active set of the solution $a^*$ to \eqref{eq:l1}, and define 
the constant $d$ as the smallest positive constant such that for any active set $\G$ visited by the LCA and any vector $x$ in $\field{R}^N$ supported on $\widetilde{\G}=\G\cup\G_*$, we have:
\begin{equation}
\prn{1-d} \eltwo{x}^2 \leq \eltwo{\Phi x}^2 \leq \prn{1+d} \eltwo{x}^2.
\label{eq:delta}
\end{equation}
Although this definition looks similar to the definition of the RIP constant, it is important to note that \eqref{eq:delta} must hold not for a general index set (as in \eqref{eq:rip})  but rather for the active sets $\widetilde{\G}$ visited by the LCA during convergence.  If such a $d$ exists, Theorem 3 of~\cite{balavoine_convergence_2012} applied to $\ell_1$-minimization becomes:
\begin{theorem}
\label{th:rate}
If $d$ defined in \eqref{eq:delta} exists and satisfies $d\leq 1$, then the LCA system \eqref{eq:dyn} converges exponentially fast with convergence speed $\prn{1-d}/\tau$, i.e. $\exists\mathcal{K}>0$, such that $\forall t\geq0$
$$\eltwo{u(t)-u^*}\leq\mathcal{K} e^{-(1-d)t/\tau}.$$
\end{theorem}
While it is difficult to characterize the path of the LCA in general, if the number of active nodes during convergence remains small, then $d$ can be related to the RIP constant and Theorem~\ref{th:rate} can be used to bound the convergence speed. This is precisely what is done in the two main results of this paper.


\section{Bounding the LCA Active Set}
\label{sec:charac}

In this section, we state our two main theorems that bound the size of the LCA active set during convergence.  
In the following, the vector $\adag\in\field{R}^N$, referred to as ``true'' underlying signal or \emph{original signal}, has $S$ non-zero coefficients supported on the set $\G_{\dagger}$, referred to as the \textit{optimal support}. This signal generates noisy measurements $y$ in $\field{R}^M$:
$$y=\Phi \adag + \epsilon = \Phi_{\G_{\dagger}} \adag + \epsilon $$
for some noise vector $\epsilon\in\field{R}^M$. Our analysis considers the general case where the measurements are corrupted by noise, but remains valid in the noise-free case, when $\epsilon = 0$.  We also define the following quantities that appear several times in the proofs of the theorems:
\begin{gather*}
 \alpha = \alpha(\delta) = (1+\delta)(1-\delta)^{-2}, \\
 C_{\delta}(p) = \alpha\prn{\eltwo{\adag} + \sqrt{1-\delta}\eltwo{\epsilon} + \lambda\sqrt{p}}.
\end{gather*}
Note that if $0<\delta<1$, then $\alpha\geq1$.


Our first result provides guarantees similar to the $S$-step property in that only the $S$ nodes that belong to the optimal support $\G_{\dagger}$ become active. Using the RIP constant $\delta$ of order $(S+1)$ of the matrix $\Phi$, the proof in Appendix~\ref{app:opt} bounds the amplitudes of the internal states $u_n(t)$ in $\Gc_{\dagger}$ to show that they remain below threshold.

\begin{theorem}
\label{th:opt}
Assume that the dictionary $\Phi$ satisfies the RIP with parameters $\prn{S+1,\delta}$ and that the support $\G(0)$ of the initial output $a(0)$ is a subset of the optimal support (i.e., $\G(0)\subset\G_{\dagger}$). If the following two conditions between the original signal $\adag$, the threshold $\lambda$, the noise $\epsilon$, the sparsity $S$ and the RIP constant $\delta$ are satisfied:
\begin{gather}
 \eltwo{\adag-a(0)} \leq C_{\delta}(S), \label{eq:icond} \\
 \prn{1-\alpha\delta\sqrt{S}}\lambda \geq \alpha\delta\prn{\eltwo{\adag} + \sqrt{1-\delta}\eltwo{\epsilon}} + \elinf{\Phi_{\Gc_{\dagger}}^T\epsilon}, \label{eq:rel1}
\end{gather}
then nodes in $\Gc_{\dagger}$ never cross threshold (i.e., $\G\subset\G_{\dagger}$).
\end{theorem}

Similarly to the analysis of some digital solvers, our second result gives conditions for the active set to be bounded by $q$, where $q$ may be larger than $S$. In contrast to the analysis for digital solvers however, these conditions do not bound the number of ``switches'' or iterations. In our case, bounding the size of the active set is enough to guarantee exponential convergence.

\begin{theorem}
\label{th:size}
Assume that the dictionary $\Phi$ satisfies the RIP with parameter $(S+q,\bar{\delta})$, for some $q\geq 0$. If the following two conditions between the original signal $\adag$, the initial state $u(0)$, the threshold $\lambda$, the noise $\epsilon$, the parameter $q$ and the RIP constant $\bar{\delta}$ are satisfied:
\begin{gather}
 \eltwo{u(0)} \leq \lambda\sqrt{q} \label{eq:initcond2} \\
\lambda \geq \dfrac{1+\bar{\delta}}{1-3\bar{\delta}}\dfrac{1}{\sqrt{q}}\bigprn{\eltwo{\adag} + \sqrt{1-\bar{\delta}}\eltwo{\epsilon}},
\label{eq:rel2}
\end{gather}
then the active set $\G$ never contains more than $q$ nodes for all time $t\geq0$ (i.e., $\abs{\G}\leq q$).
\end{theorem}

The proof in Appendix~\ref{app:size} bounds the energy of the $q$ biggest nodes $u_n(t)$. The simulations in Section~\ref{sec:sim} show that useful values for $q$ are typically small multiples of $S$. Conditions \eqref{eq:icond}, \eqref{eq:rel1} and \eqref{eq:rel2} involve complex relationships between the various parameters. We analyze their implication on the RIP constant below, and look at the resulting number of measurements for CS random matrices in Section~\ref{ssec:discussion}.

First, condition \eqref{eq:icond} of Theorem~\ref{th:opt} constrains the starting point to be reasonably close to the optimum $\adag$. When the system starts at rest, $u(0)=0$ and condition \eqref{eq:icond} becomes:
$$\eltwo{\adag}\leq \alpha\bigprn{\eltwo{\adag}+\sqrt{1-\delta}\eltwo{\epsilon}+\lambda\sqrt{S}},$$
which always holds since $\alpha\geq1$. Similarly, if the system starts at rest, condition~\eqref{eq:initcond2} obviously holds.

To analyze the requirement of the theorems on the RIP more easily, we assume without loss of generality that $\eltwo{\adag}=1$. For now, we also assume that there is no noise, $\epsilon=0$; the noise level is addressed in Section~\ref{ssec:discussion}.
From \eqref{eq:opt}, the solution $a^*$ is a thresholded version of the original signal $\adag$:
$$a^*_{\G_*} = \adag_{\G_*} - \lambda \prn{\Phi_{\G_*}^T\Phi_{\G_*}}^{-1} z_{\G_*}.$$
Even though the two theorems require the threshold $\lambda$ to be sufficiently large, it cannot be taken too large or the outputs simply remain zero. However, if some nodes in $\G_{\dagger}$ have small amplitudes, they do not contribute much to the signal energy and setting them to zero in $a^*$ may be acceptable.

It is instructive then to look at the scenario where all the non-zero entries in $\adag$ have the same magnitude. In this case, they contribute equally to the signal energy and the threshold should be low enough to recover them all. If $\eltwo{\adag} = 1$, each non-zero element of $\adag$ is equal to $\pm 1/\sqrt{S}$. Thus, the threshold should be below $1/\sqrt{S}$. We take $\lambda = {r}/{\sqrt{S}}$, for some $0<r<1$.
Rearranging the terms in \eqref{eq:rel1} yields the following condition on the RIP constant in Theorem~\ref{th:opt}:
\begin{equation}
 \delta \leq \dfrac{r}{\prn{1+r}\alpha\sqrt{S}}.
 \label{eq:deltacond1}
\end{equation}
This shows that, for the active set to remain a subset of the optimal support, the RIP constant needs to scale with $1/\sqrt{S}$. 
%
%

Taking $q=\beta S$ in Theorem~\ref{th:size}, for a constant $\beta$, \eqref{eq:rel2} becomes
\begin{equation}
 \bar{\delta} \leq \dfrac{r\sqrt{\beta} - 1}{3r\sqrt{\beta} +1}.
 \label{eq:deltacond2}
\end{equation}
This shows that for the active set to contain less than $q$ nodes, the RIP constant needs only to be bounded by a small constant that does not depend on $S$ anymore, which is more favorable than condition~\eqref{eq:deltacond1}.

\section{Application to Compressed Sensing matrices}
\label{ssec:discussion}

Theorems~\ref{th:opt} and~\ref{th:size} are deterministic. However, to get a better feel for their implications on the various parameters, the results can be interpreted in terms of a known expression for the RIP constant established for random matrices. In particular, Theorem 5.65 in \cite{vershynin_introduction_2011} states that if $\Phi$ is an $M\times N$ random matrix, whose columns $\Phi_n$ are independent subgaussian random vectors in $\field{R}^M$ with $\eltwo{\Phi_n}=1$, then for any sparsity level $1\leq S\leq N$ and any $\delta\in\prn{0,1}$, the matrix $\Phi$ satisfies the RIP with parameters $\prn{S,\delta}$ with high probability provided
\begin{equation}
\delta \sim \sqrt{\dfrac{S \log (N/S)}{M}},
\label{eq:ripest}
\end{equation}
where $\sim$ indicates that two quantities are equal up to a constant.
Some examples are random matrices with independent and identically distributed Bernoulli columns with unit norm, or columns drawn independently and uniformly at random from the unit sphere. This result implies that typical CS matrices have a number of measurement $M$ on the order of $ \mathcal{O}\prn{{\delta^{-2}}S\log\prn{{N}/{S}}}$.

\subsection{Theorem~\ref{th:opt} with CS matrices}

\subsubsection*{\textbf{Measurements}}
Plugging the estimate \eqref{eq:ripest} for $\delta$ in \eqref{eq:deltacond1} yields:
$$ \sqrt{M} \gtrsim S\sqrt{\log{N/S}} ~\dfrac{(1+r)\alpha}{r}, $$
where the notation $\gtrsim$ means greater up to a constant factor. When $S \ll M$, $\delta$ is small and $\alpha \sim 1$. As a reference, for an RIP of $\delta\leq 1/2$, $\alpha\leq6$, and for $\delta\leq0.1$, $\alpha\leq1.358$. This shows that the number $M$ of measurements for a random matrix $\Phi$ must be on the order of $\mathcal{O}(S^2 \log{(N/S)})$.

This result strongly resembles the condition for the Homotopy algorithm to satisfy the $S$-step property~\cite{donoho_fast_2008}, which requires that $S \leq \prn{1+\mu^{-1}}/{2}$ and leads to the same number of measurements. For $M\!\sim \mathcal{O}(S^2 \log N)$, the Homotopy algorithm on the parameter $\lambda$ behaves like a pursuit algorithm, where nodes are added to the active support and the solution evolves in a piecewise-linear manner. Likewise, the LCA solution evolves according to a continuous switched linear system and nodes are added to the support until the solution is reached. Both results ensure that only nodes present in the original signal enter the active set.
OMP was also shown to recover an $S$-sparse signal in exactly $S$ steps provided that $\Phi$ satisfies the RIP with $\sqrt{S} \leq {1}/\prn{3\delta_{S+1}}$, which also leads to $\mathcal{O}(S^2\log N)$ measurements~\cite{davenport_analysis_2010}. Consequently, despite the continuous-time nature of the LCA trajectories, comparable bounds on the RIP constant emerge from our study.

\subsubsection*{\textbf{Noise level}}
When $\epsilon$ is a Gaussian white noise whose entries have variance $\sigma^2$, the terms due to the noise in \eqref{eq:rel1} are $\eltwo{\epsilon} \sim \sqrt{M}\sigma$ and $\elinf{\Phi^T\epsilon} \sim \sqrt{\log N}\sigma$ with high probability. Taking these terms into account in \eqref{eq:rel1} does not change the bound on $\delta$ in \eqref{eq:deltacond1} by more than a constant if
$$\alpha\delta\sqrt{1-\delta}\eltwo{\epsilon}+\eltwo{\Phi^T_{\Gc_{\dagger}}\epsilon} = \kappa \alpha\delta,$$
for some constant $\kappa >0$. Using the estimate \eqref{eq:ripest}, along with $S\ll N$, $M\sim S^2\log(N/S)$, $\alpha \sim 1$, $\alpha\sqrt{1-\delta}\sim 1$, and reorganizing the terms yield a noise variance of:
\begin{align*}
\sigma  \sim \dfrac{\alpha\delta\kappa}{\alpha\delta\sqrt{1-\delta}\sqrt{M}+\sqrt{\log N}} & \sim \dfrac{\kappa \sqrt{\frac{S\log(N/S)}{M}}}{\sqrt{S\log(N/S)} + \sqrt{\log N}} \\
  \sim \dfrac{\kappa}{1+\sqrt{\frac{\log N}{S\log (N/S)}}}\dfrac{1}{\sqrt{M}} & \sim \dfrac{\kappa}{1+\frac{1}{\sqrt{S}}}\dfrac{1}{\sqrt{M}}.
\end{align*}
Thus, the total energy allowed in the noise vector is on the order of $\eltwo{\epsilon} \sim \prn{1+1/\sqrt{S}}^{-1}$, which is approximately on the same order as the energy of the signal.

This result can be improved upon. Theorem~\ref{th:opt} is stated for any fixed noise vector $\epsilon$. In the case where the noise $\epsilon$ is assumed to be a Gaussian random vector, the proof of Lemma~\ref{lem:bounddist} in Appendix~\ref{app:lemma} hints that the bound used for $\eltwo{a^{\infty}-\adag}$ can be improved. The essential step consists in bounding $\eltwo{(\PP)^{-1}\Phi_{\G}^T\epsilon}$. It is an easy calculation that
$$ \E{\eltwo{(\Phi_\Gamma^T\Phi_\Gamma)^{-1}\Phi_\Gamma^T \epsilon}^2} = 
\sigma^2 \operatorname{Trace}\prn{(\Phi_\Gamma^T\Phi_\Gamma)^{-1}}
\leq \frac{S\sigma^2}{1-\delta}.$$
Moreover, standard tail inequalities~\cite{vershynin_introduction_2011} show that this random variable concentrates around its mean. Thus, when the noise is Gaussian, $\sqrt{1-\delta}\eltwo{\epsilon}$ can be replaced by $ \sqrt{S}\sigma$ with high probability in \eqref{eq:rel1}. Going over the equations in the previous paragraph, we obtain a noise variance of the form:
$$\sigma \sim {\dfrac{\alpha\delta\kappa}{\alpha\delta\sqrt{S}+\sqrt{\log N}}} \sim \dfrac{\kappa}{1+\sqrt{\log N}}\dfrac{1}{\sqrt{S}}.$$
The total energy allowed in the noise vector becomes $\eltwo{\epsilon} \sim \left.\sqrt{M/S}\middle/\prn{1+\sqrt{\log N}}\right.$, which increases with the number of measurements $M$.

\subsection{Theorem~\ref{th:size} with CS matrices}

\subsubsection*{\textbf{Measurements}}
For random matrices, using the estimate \eqref{eq:ripest} of the RIP constant $\bar{\delta}$ of order $S+q = (1+\beta)S$ in \eqref{eq:deltacond2} yields
$$ \sqrt{M} \gtrsim \sqrt{(1+\beta)S\log\prn{\dfrac{N}{(1+\beta)S}}} ~\dfrac{3r\sqrt{\beta}+1}{r\sqrt{\beta}-1}.$$
For $\beta$ constant, this yields a number of measurements on the order of $\mathcal{O}\prn{S\log{(N/S)}}$. 
For reference, using $\beta = 30$ and $r=0.8$ in \eqref{eq:deltacond2} yields $\delta_{31S}\leq0.23$. In comparison, OMP has been shown to converge for $\delta_{31S}\leq 1/3$ \cite{zhang_sparse_2011}. The result for ROMP has a slightly worse form with $\delta_{8S}\leq 0.01/\sqrt{\log S}$ \cite{needell_signal_2010}. Finally, CoSaMP was shown to converge for $\delta_{4S}\leq0.1$ \cite{needell_cosamp:_2008}. For all those algorithms, the RIP constants reported lead to the same order of measurements $\mathcal{O}\prn{S \log N}$. This is another interesting parallel between the LCA and its digital equivalents. Letting more than the optimal $S$ nodes enter the active support still yields good convergence results, while giving better scaling on the RIP constant and number of measurements. Contrary to the digital solvers however, the conditions for the LCA are only necessary to guarantee a bound on the exponential speed of convergence. As stated before, the algorithm is guaranteed to converge to the solution to \eqref{eq:l1} without any requirements on the RIP constant. Moreover, the error achieved by the LCA is linked to the performance guarantees associated with $\ell_1$-minimization, as discussed in Section~\ref{sec:algo}.

\subsubsection*{\textbf{Noise level}}
The influence of the noise appears clearly in our results. The noise term in \eqref{eq:rel2} does not affect the bound on $\bar{\delta}$ in \eqref{eq:deltacond2} by more than a constant if
$$\sqrt{1-\bar{\delta}}\eltwo{\epsilon} = \bar{\kappa},$$
for some $\bar{\kappa}>0$. Assuming again that $\epsilon$ is a Gaussian white noise, whose entries have variance $\sigma^2$, and that $\eltwo{\adag}=1$ yields a noise variance of
$$\sigma \sim \dfrac{\bar{\kappa}}{\sqrt{1-\bar{\delta}}}\dfrac{1}{\sqrt{M}} \sim \dfrac{1}{\sqrt{M}}.$$
As a consequence, the total energy $\eltwo{\epsilon}$ allowed in the noise vector is $\mathcal{O}\prn{1}$, which is the same order as the energy of the signal. Here again, assuming that the noise is Gaussian in the proof of the theorem itself leads to a sharper bound. Using the same concentration argument as previously, we can replace $\sqrt{1-\delta}\eltwo{\epsilon}$ by $\sqrt{q}\sigma$ with high probability in \eqref{eq:rel2}. This yields a new noise variance of the form $\sigma \sim {\kappa}/{\sqrt{q}}$ and the energy in the noise vector becomes $\eltwo{\epsilon} \sim \sqrt{M/q}$, which can again increase with the number of measurements.

\subsection{Decreasing threshold}
\label{ssec:threshdecay}

Interestingly, the proofs of the two theorems hint that the results possibly still hold when the threshold is not constant, but rather exponentially decreasing. 
In the proof of Theorem~\ref{th:opt}, the lower bound on the threshold $\lambda$ depends on the quantity $\eltwo{a(t)-\adag}$, and in the proof of Theorem~\ref{th:size}, it depends on $\eltwo{u(t)-u^*}$. If the system is exponentially convergent, both quantities should decrease exponentially fast with time. Thus, the threshold $\lambda$ could be decreased according to an exponential decay while still satisfying the inequalities in the two theorems. Decreasing the threshold would allow the system to potentially recover more nodes from $\adag$, while keeping the size of the active set bounded, as well as yielding faster convergence. This is confirmed in simulation (see Section~\ref{sec:sim}). Interestingly, similar observations have been made for digital solvers (e.g. in~\cite{hale_fixed-point_2008}, the threshold is decreased according to a geometric progression to speed up recovery). However, there has been no analytic justification for the observed increase in speed or for how to choose the decay rate. In our case, even if the proof suggests the potential advantage of decreasing the threshold according to an exponential decay, the additional dynamics on the threshold would drastically change the nature of the analysis, starting with the proof of convergence in~\cite{balavoine_convergence_2012}.

\subsection{Estimate of the Convergence speed}
\label{ssec:speed}

The ultimate goal of this paper is to obtain an estimate of the speed of convergence of the LCA in the context of CS recovery. In Theorem~\ref{th:opt}, we showed that under some conditions, the active sets visited during convergence may never contain more than the $S$ optimal nodes. This result was generalized in Theorem~\ref{th:size} to allowing no more than $q$ nodes to become active, where $q$ is typically a small multiple of $S$. Such guarantees allow us to approximate $d$ in Theorem~\ref{th:rate} by the RIP constant of $\Phi$ of order $S$ or $q$. Thus, for random matrices of interest, $d$ is approximately $\sqrt{S\log (N/S) /M}$. The convergence being exponential, this leads to an estimate for the convergence time of the LCA of $\mathcal{O}\prn{\dfrac{\tau}{1-\sqrt{S\log (N/S)/M}}}$, where $\tau$ is the time constant of the physical solver.

For informational purposes, the digital solvers Homotopy, OMP, ROMP and CoSamp have been proven to have running times on the order of $\mathcal{O}(SMN)$ flops when the number of iterations is finite \cite{donoho_fast_2008,zhang_sparse_2011,needell_signal_2010,needell_cosamp:_2008}. This estimate can in general be reduced if a fast multiply for $\Phi$ and $\Phi^T$ is available. It is important to keep in mind that the time constant $\tau$ for the LCA has the potential to be much smaller than the time to perform a single matrix multiply for a digital solver \cite{shapero_low_2012}. As a consequence, the scaling properties of the LCA seem more favorable for large problems.

\begin{figure}
  \centering
	  \includegraphics[width=2.5in]{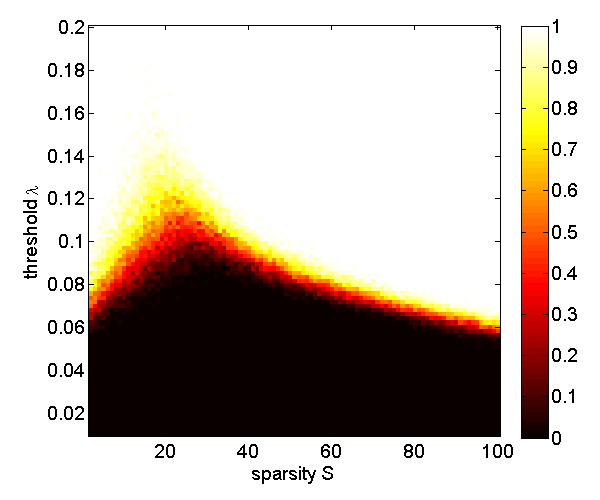}
	  \caption{Percentage of the trials where no more than the $S$ nodes from the optimal support $\G_{\dagger}$ become active during convergence. The value 1 means that 100\% of the trials satisfied this condition.}
	  \label{fig:opt}
\end{figure}


\section{Simulations}
\label{sec:sim}

In this section, we provide simulations that illustrate the previous theoretical results\footnote{Matlab code running the experiments in this section can be downloaded at \url{http://users.ece.gatech.edu/~abalavoine3/code/LCA_CS_exp.zip}}. As an example, we use a sparse vector $\adag$ of length $N=400$ whose non-zero entries are generated by randomly selecting $S=5$ indices and setting their amplitudes so that $\eltwo{\adag}=1$. Then, we take $M=200$ measurements by generating a Gaussian random matrix $\Phi$ of size $200\times 400$, with entries drawn independently from a normal distribution and columns normalized to have unit norm. We also add a Gaussian white noise with standard deviation $\sigma=0.025$ to the measurement so that $y=\Phi \adag + \epsilon$. All the results of this section are obtained by simulating the ODEs \eqref{eq:dyn} on a digital computer. The algorithm is always started at rest with $u(0)=0$.

\subsection{Effect of the threshold on the size of the active set}
\begin{figure}
\centering
	  \includegraphics[width=2.5in]{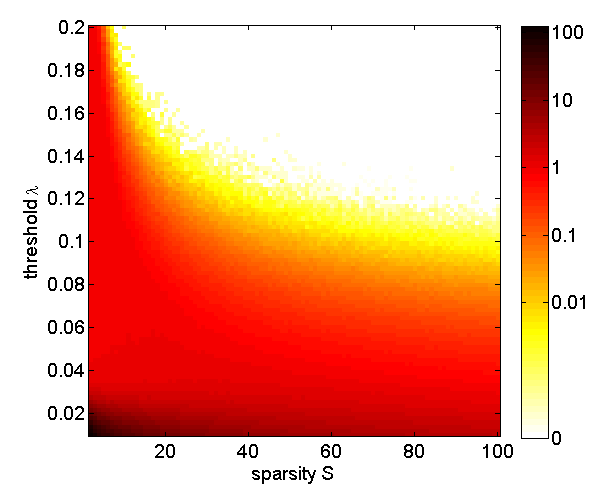}
	  \caption{Ratio of the maximum number of active elements $q$ during convergence over the sparsity level $S$. For instance, a value of $10$ in the color bar means that the biggest active set during convergence contains $10S$ active elements.}
	  \label{fig:ratio}
\end{figure}

We first explore how the value of the threshold $\lambda$ affects the size of the active set during convergence, i.e. the maximum number of nodes that become active while the system is evolving. In \figurename~\ref{fig:opt} and \ref{fig:ratio}, we vary the value of the threshold $\lambda$ and the sparsity level $S$. For each point on the figures, we simulate 100 random draws of a sparse vector $\adag$ and a measurement matrix $\Phi$ and assume that no noise is present. In \figurename~\ref{fig:opt}, we look at the percentage of the 100 trials where only nodes that are part of the optimal support $\G_{\dagger}$ become active.  For large $S$ (approximately $S>28$), the transition phase for $\lambda$ follows a curve that looks like $1/\sqrt{S}$. For small $S$, the behavior appears qualitatively different. Both follow the general prediction from \eqref{eq:rel1}:
$$\lambda \gtrsim \dfrac{\alpha\delta}{1-\alpha\delta\sqrt{S}}.$$
In \figurename~\ref{fig:ratio}, the color coding represents the ratio of the maximum number of active elements $q$ during convergence over the sparsity level $S$. The phase transition on this plot follows a $1/\sqrt{S}$ behavior. Moreover, for most values of the threshold $\lambda$ and sparsity level $S$, a relatively few number of elements become active during convergence ($q$ is mostly contained between $1S$ and $10S$).  The results in \figurename~\ref{fig:opt} and \figurename~\ref{fig:ratio} confirm the qualitative behavior of the bounds derived in Theorems~\ref{th:opt} and \ref{th:size}.

\begin{figure}
\centering
\includegraphics[width=3.4in]{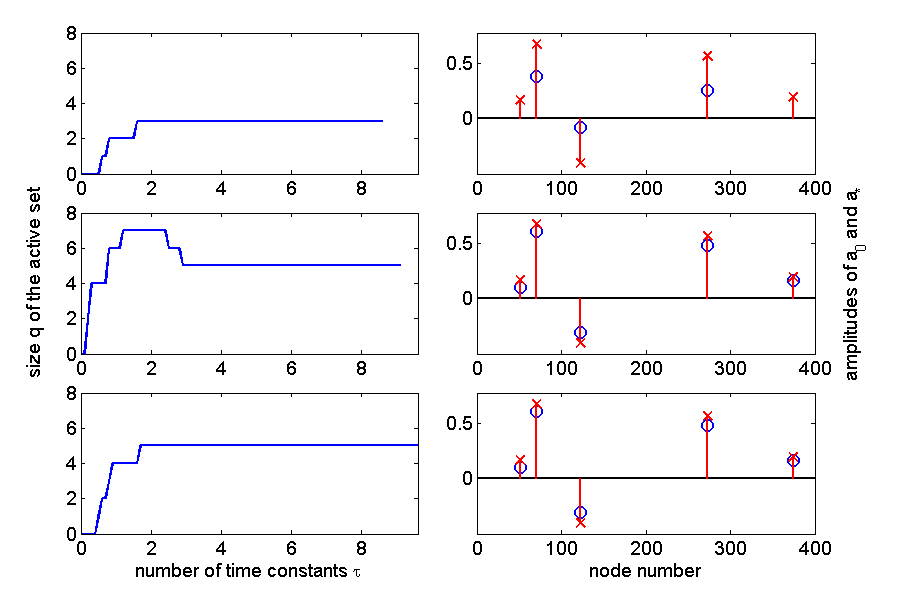}
\caption{This figure shows the number of active nodes (left column) and the fixed point $a_*$ reached by the LCA (right column), for different choices of the threshold.  The red crosses represent the original signal $\adag$ and the blue rounds are the solutions $a_*$. A fixed threshold $\lambda=0.3$ was used in the first row, $\lambda=0.08$ in the second row, and the threshold was decreased from $0.3$ to $0.08$ according to an exponential decay in the third row. }
\label{fig:decayl}
\end{figure}

\begin{figure*}
\centering
	\subfigure[Effect of the signal length $N$.]{
  \includegraphics[width=2.3in]{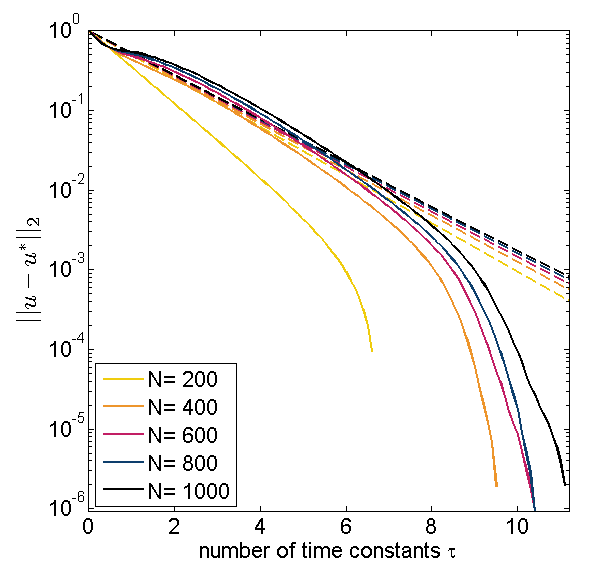}
	\label{fig:varyN} }
	\subfigure[Effect of the sparsity level $S$.]{
  \includegraphics[width=2.3in]{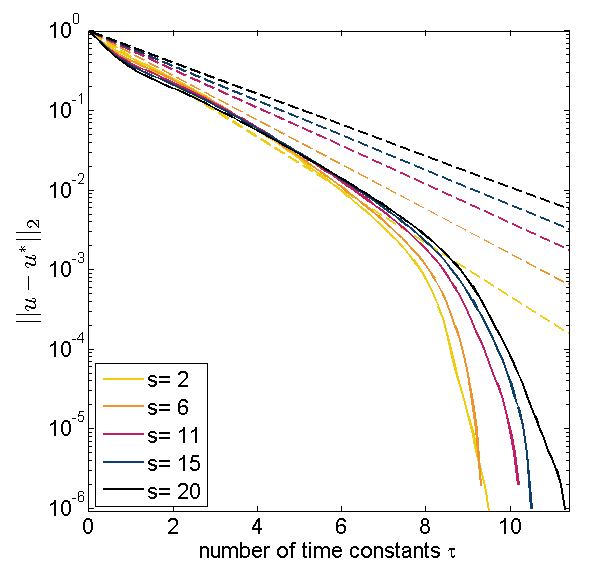}
	\label{fig:varyS} }
	\subfigure[Effect of the number of measurements $M$.]{
  \includegraphics[width=2.3in]{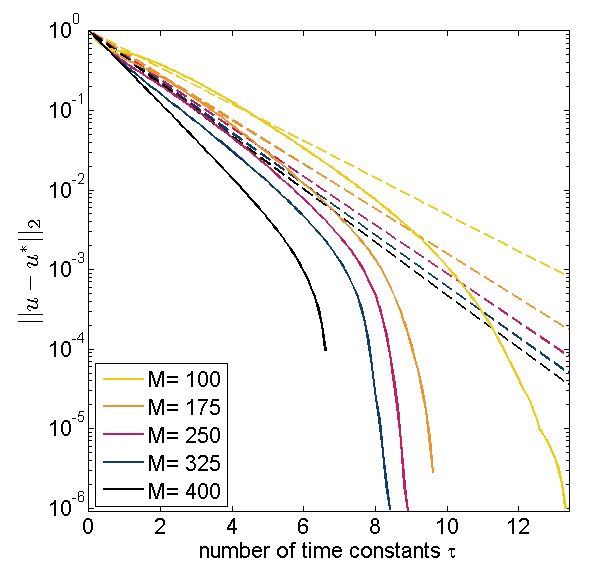}
	\label{fig:varyM} }
	\subfigure[Effect of the threshold $\lambda$.]{
	\includegraphics[width=2.3in]{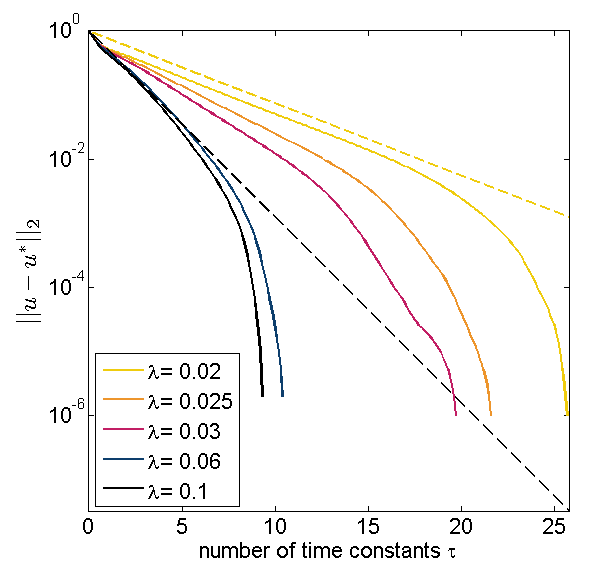}
	\label{fig:varyl} }
	\caption{Evolution of the experimental mean-squared error $\eltwo{u(t)-u^*}$ (plain line) and theoretical decay (crossed line) as problem parameters are varied.}
	\label{fig:vary1}
\end{figure*}

\subsection{Decreasing the threshold during convergence}

As mentioned in Section~\ref{ssec:threshdecay}, the proofs of Theorems~\ref{th:opt} and \ref{th:size} suggest that the active set remains bounded even when the threshold is decreased according to an exponential decay, while yielding faster convergence. This fact is confirmed in practice. To illustrate this, we first ran the LCA with a high threshold value of $\lambda=0.3$. As shown in the first row of \figurename~\ref{fig:decayl}, the active set never contains more than three nodes that are part of the optimal support. The final solution is missing two nodes from the original signal $\adag$. In the second row, $\lambda$ is fixed to a low value of $0.08$. The final solution recovers all the nodes from $\adag$. However, the biggest active set now contains $q=7$ nodes and the convergence is slower. Finally, in the last row, the threshold is started at $0.3$ and decreased to the value $0.08$ according to an exponential decay. As expected, the final solution is the same as the one in row 2. However, in this case the active set never contains more than the five nodes from the optimal support. Moreover, the system converges faster, in less than $2\tau$ compared to $3\tau$ in row 2.

\subsection{Estimate of the convergence speed}

Finally, we would like to know how well the quantity
\begin{equation}
e^{-(1-\delta)t/\tau},
\label{eq:thdecay}
\end{equation}
predicted by Theorem~\ref{th:rate} bounds the convergence rate of the solver, represented by the mean-squared error between the nodes at time $t$ and the final solution $u^*$:
$$\eltwo{u(t)-u^*}.$$
In \figurename~\ref{fig:vary1}, this quantity is normalized to start at $1$ so as to compare it with \eqref{eq:thdecay}. When they are not varying, we fix the threshold $\lambda=0.1$, the number of measurements $M=200$, the sparsity $S=5$, and the signal length $N=400$. For each experimental curve (solid lines), the mean-squared error is averaged over 100 trials. We also plot the theoretical decay (dashed lines) using the expression \eqref{eq:thdecay} with $\delta = \sqrt{{S\log (N/S)}/{M}}$.
As expected, the theoretical curves approximate the decay of the experimental mean-squared error. These upper bounds are not strict in practice since only an estimate for the RIP constant $\delta$ can be used. However, this is enough to check that the experimental curves qualitatively follow the theoretical predictions as the parameters $N$, $M$ or $S$ are varied in \figurename~\ref{fig:varyN}, \ref{fig:varyM} and \ref{fig:varyS} respectively.

In \figurename~\ref{fig:varyl}, we explore the effect of the threshold $\lambda$ on the experimental decay. For values of $\lambda$ bigger than $0.06$, the bound \eqref{eq:thdecay} with $\delta = \sqrt{{S\log (N/S)}/{M}}$ (dark blue dashed line) is valid, even though more than $S$ nodes may become active (for $\lambda=0.06$, the average over 100 trials for the maximum size for the active is $q=23=4.6 S$). As $\lambda$ becomes smaller, more nodes are able to enter the active set. The estimate for $\delta$ must be changed to $\delta = \sqrt{q\log(N/q)/M}$, where $q$ is the maximum number of active elements during convergence. Using $\delta = \sqrt{5S \log (N/S) /M}$ (yellow dashed line) provides an upper bound even for very small values of the threshold $\lambda$, where much more than $5S$ nodes become active during convergence (the average over 100 trials for the maximum size of the active set for $\lambda=0.02$ is $180=36 S$).

\section{Conclusions}
\label{sec:ccl}
In this paper, we studied a dynamical system for solving $\ell_1$-minimization problems in the context of CS signal recovery. In this specific problem setting we are able to give strong guarantees about the path followed by the system's internal state variables during convergence. Indeed, our results show that in typical CS situations, the path followed by the LCA is close to optimal, with only a few nodes entering the active set during convergence.  These results can then be used to make strong guarantees on the exponential convergence speed of the system, and the quantitative results generally agree qualitatively with our simulation results.  Interestingly, despite the LCA being a completely different computing architecture than traditional algorithms being run on a digital computer, the conditions of our results directly parallel the established guarantees for several digital algorithms.   As with any signal processing system, such performance guarantees are important to establish before investing significant resources in system development and deploying the system in an application.  The strong performance guarantees of this paper lead us to conclude that the LCA, if implemented in a large-scale analog circuit, could lead to substantial improvements in the time required for CS signal recovery in many problems of interest.  


\appendices


\section{Nodes Dynamics}
\label{app:ode}

The LCA is a type of switched linear system \cite{decarlo_perspectives_2000}, where the dynamics are a linear ODE that changes every time a node crosses threshold (i.e., moves into or out of the active set).  Between switching times, the active set $\G$ is fixed, $\dot{a}_{\G}(t) = \dot{u}_{\G}(t)$, and the ODE \eqref{eq:dyn} can be partially decoupled:
\begin{align} 
 \dot{a}_{\G}(t) & = -\PP a_{\G}(t) + \Phi_{\G}^Ty -\lambda z_{\G}(t),
\label{eq:active} \\
 \dot{u}_{\Gc}(t) & = -u_{\Gc}(t)-\PcP a_{\G}(t) + \Phi_{\Gc}^Ty .
\label{eq:inactive}
\end{align}

The following is a very brief overview of some fundamental results in linear ODEs, which we use freely in our proofs. Let $x(t)$ be a function from $\field{R}^+$ to $\field{R}^N$, $A$ a symmetric matrix in $\field{R}^{N\times N}$ and $b$ a vector in $\field{R}^N$. The solution to 
\begin{equation}
\dot{x}(t) = Ax(t) + b.
\label{eq:ode1}
\end{equation}
with initial condition $x(t_k)=x^{t_k}$ is
$$x(t) = e^{A (t-t_k)}x^{t_k} +\prn{I-e^{A(t-t_k)}}A^{-1}b.$$

The above expression $\prn{I-e^{At}}A^{-1}$ is always well-defined even when the matrix $A$ is singular. To see this, first diagonalize the matrix as $A = P\Lambda P^{-1}$, where $\Lambda$ is a diagonal matrix with diagonal elements $\lambda_i$: $\Lambda=\mathrm{diag}\prn{\lambda_1,\ldots,\lambda_n}$. Plugging this in the above expression yields:
\begin{align*}
\left( I \right. & \left. -e^{At} \right) A^{-1} = P\prn{I-e^{\Lambda t}}\Lambda^{-1} P^{-1}\\
\ & = P\mathrm{diag}\Bigl( \prn{1-e^{\lambda_1 t}}\lambda_1^{-1} ,\ldots, \prn{1-e^{\lambda_n t}}\lambda_n^{-1} \Bigr) P^{-1}
\end{align*}
To see that diagonal elements are still defined when $\lambda_i =0$, we first take a Taylor expansion when $\lambda_i$ goes to zero:
$$\lambda_i^{-1}\prn{1-e^{\lambda_i t}} = \lambda_i^{-1} \Bigl(-\lambda_i t + o(\lambda_i^2)\Bigr) = -t + o(\lambda_i).$$
By continuity, we get that $\prn{1-e^{\lambda_i t}}\lambda_i^{-1} = -t$ when $\lambda_i=0$. Thus, the matrix $\prn{I-e^{At}}A^{-1}$ is well defined.

In the case where $b$ varies with time, the solution to \eqref{eq:ode1} with initial condition $x(t_k)=x^{t_k}$ is:
\begin{equation}
 x(t) = e^{A(t-t_k)} x^{t_k} + e^{A t}\dint{t_k}{t}{e^{-A\nu} b(\nu) d\nu}.
\label{eq:odesol2}
\end{equation}

Applying these results, the solution to \eqref{eq:active} on the active set $\G$ between switching times $t_k$ and $t_{k+1}$ is given by
\begin{equation}
a_{\G}(t) = e^{-A(t-t_k)}a_{\G}^{t_k} + \prn{I- e^{-A(t-t_k)}} A^{-1}\prn{\Phi_{\G}^T y - \lambda z_{\G}},
\label{eq:actsol}
\end{equation}
where $A=\PP$ and $a_{\G}^{t_k} = a_{\G}(t_k)$.
In the case where $\PP$ is non-singular, the term $a_{\G}^{\infty} = A^{-1}\prn{\Phi_{\G}^T y - \lambda z_{\G}}$ can be interpreted as the steady state of \eqref{eq:active} if the active set and sign vector $z_{\G}$ remain unchanged until convergence. The points $a_{\G}^{\infty}$ play a key role in our proofs (see Lemma~\ref{lem:bounddist}).

Finally, the solution to the linear ODE \eqref{eq:inactive} on the inactive set $\Gc$ between switching times $t_k$ and $t_{k+1}$ is given by:
\begin{equation}
u_{\Gc}(t) = e^{-(t-t_k)}u_{\Gc}^{t_k} + e^{-t}\dint{t_k}{t}{e^{\nu} \rho_{\Gc}(\nu)d{\nu}},
\label{eq:inactsol1}
\end{equation}
where $\rho_{\Gc}(\nu) = \Phi_{\Gc}^T\prn{y-\Phi_{\G}a_{\G}(\nu)}$ and $u_{\Gc}^{t_k} = u_{\Gc}(t_k)$.
Letting $t$ go to infinity in equations \eqref{eq:actsol} and \eqref{eq:inactsol1}, the fixed point $a^*$ supported on the final active set $\G_*$ must satisfy:
\begin{align*}
 & a^*_{\G_*} = \prn{\Phi_{\G_*}^T\Phi_{\G_*}}^{-1}\prn{\Phi_{\G_*}^T y - \lambda z_{\G_*}} \\
 & u^*_{\Gc_*} = \Phi_{\Gc_*}^T\prn{y-\Phi_{\G_*}a^*_{\G_*}}.
\end{align*}
Since a node $j$ is in the inactive set $\Gc_*$ if and only if $\abs{u_j}\leq\lambda$, the two equations above translate immediately to:
\begin{equation}
\begin{split}
 & a^*_{\G_*} = \prn{\Phi_{\G_*}^T\Phi_{\G_*}}^{-1}\prn{\Phi_{\G_*}^T y - \lambda z_{\G_*}} \\
 & \elinf{\Phi_{\Gc_*}^T\prn{y-\Phi_{\G_*}a^*_{\G_*}}}\leq\lambda,
\end{split}
\label{eq:opt}
\end{equation}
which are the two well-known optimality conditions for $a^*$ to be the solution to \eqref{eq:l1} \cite{fuchs_sparse_2004}. 


\section{Lemmas}
\label{app:lemma}

The proofs of Theorems~\ref{th:opt} and \ref{th:size} make use of the following two lemmas. The first lemma bounds the $\ell_2$-distance between the points $a^{\infty}_{\G}$ and the true signal $\adag$.

\begin{lemma}
\label{lem:bounddist}
Let $a^{\infty}$ be a vector supported on a set $\G$ that contains less than $p$ indices and that satisfies:
$$\PP a^{\infty} = \Phi_{\G}^Ty - \lambda z_{\G},$$
where $z_{\G} = \sign{a_{\G}^{\infty}}$. Let $R=\abs{\G\cup\G_{\dagger}}$ be the number of elements in the support of $\prn{a^{\infty}-\adag}$. If $\Phi$ satisfies the RIP with parameters $(R, \delta)$, then the following holds:
$$\eltwo{a^{\infty}-\adag} \leq \underbrace{(1-\delta)^{-1}\prn{\eltwo{\adag} + \sqrt{1-\delta} \eltwo{\epsilon} + \lambda\sqrt{p}}}_{= (1-\delta)(1+\delta)^{-1}C_{\delta}(p)}.$$
\end{lemma}

\begin{IEEEproof}
We start by noting that since $\Phi$ satisfies the RIP and $p\leq R$, then $\norm{\prn{\PP}^{-1}}\leq(1-\delta)^{-1}$, $\norm{\Phi_{\G}^T\Phi_{\G_{\dagger}\cap\Gc}}\leq\delta$ as a submatrix of $\Phi^T\Phi - I$ with at most $S\leq R$ columns, and $\norm{\prn{\PP}^{-1}\Phi_{\G}^T}^2 \leq (1-\delta)^{-1}$ \cite[Prop. 3.1, 3.2]{needell_cosamp:_2008}.\\
Splitting $\adag$ into its component on $\G$ and $\Gc$, we get that:
\begin{align*}
& \adag_{\G} = \prn{\PP}^{-1}\PP \adag_{\G}, \\
& \Phi\prn{\adag - \adag_{\G}} = \Phi_{\Gc}\adag_{\Gc}.
\end{align*}
We use these facts to finish the proof:
\begin{align*}
 & \eltwo{a^{\infty} - \adag} = \eltwo{\prn{\PP}^{-1}\prn{\Phi_{\G}^Ty-\lambda z_{\G}} - \adag}  \\
 & ~~ = \eltwo{ \prn{\PP}^{-1}\prn{\Phi_{\G}^T\prn{\Phi \adag+\epsilon} - \lambda z_{\G} } -\adag_{\G} - \adag_{\Gc}} \\
 & ~~ = \left\| \prn{\PP}^{-1} \Phi_{\G}^T\Phi_{\Gc} \adag_{\Gc} + \adag_{\G} + \prn{\PP}^{-1}\Phi_{\G}^T\epsilon  \right. \\
 & \qquad \qquad \left.  - \lambda \prn{\PP}^{-1} z_{\G} - \adag_{\G} - \adag_{\Gc} \right\|_2 \\
 & ~~ \leq \norm{\prn{\PP}^{-1}} \norm{\Phi_{\G}^T\Phi_{\G_{\dagger}\cap\Gc}} \eltwo{\adag_{\Gc}} + \eltwo{\adag_{\Gc}} \\
 & \qquad + \norm{\prn{\PP}^{-1}\Phi_{\G}^T}\eltwo{\epsilon} + \lambda\norm{\prn{\PP}^{-1}}\eltwo{z_{\G}}  \\
 & ~~ \leq (1-\delta)^{-1}\delta \eltwo{\adag_{\Gc}} +  \eltwo{\adag_{\Gc}}  \\
 & \qquad \qquad + \sqrt{1-\delta}^{-1}\eltwo{\epsilon}  + \lambda(1-\delta)^{-1}\sqrt{p} \\
 & ~~ \leq (1-\delta)^{-1} \prn{\eltwo{\adag} + \sqrt{1-\delta}\eltwo{\epsilon} + \lambda\sqrt{p}}. \tag*{\IEEEQED}
\end{align*}
\let\IEEEQED\relax
\end{IEEEproof}

The second result states that the $\ell_2$-distance of the output $a(t)$ to the true signal $\adag$ remains bounded for all time $t\geq0$.

\begin{lemma}
\label{lem:timebound}
Assume that at switching time $t_k$, the current active set $\Gk$ contains less than $p$ indices, that $\Phi$ satisfies the RIP with parameters $(R_k,\delta)$, where $R_k=\abs{\Gk\cup\G_{\dagger}}$, and that
$$\eltwo{a(t_k)-\adag} \leq C_{\delta}(p).$$
Then, for all $t\in\Prn{t_k,t_{k+1}}$,
$$\eltwo{a(t)-\adag}\leq C_{\delta}(p).$$
\end{lemma}

\begin{IEEEproof}
Define $a_{\Gk}^{\infty}$ by $\PPk a_{\Gk}^{\infty} = \Phi_{\Gk}^Ty-\lambda z_{\Gk}$. Applying Lemma~\ref{lem:bounddist}, we obtain that
$$\eltwo{a_{\Gk}^{\infty}-\adag}\leq (1-\delta)(1+\delta)^{-1}C_{\delta}(p).$$
Using the dynamics in \eqref{eq:actsol}, we have that for $t\in[t_k,t_{k+1})$:
\begin{align*}
 &  \eltwo{a(t)-\adag} = \eltwo{a_{\Gk}(t)-\adag} \\
 & \qquad = \left\| e^{-A\prn{t-t_k}}a_{\Gk}(t_k) \right.  + \left. \prn{I-e^{-A\prn{t-t_k}}}a_{\Gk}^{\infty} -\adag \right\|_2 \\
 & \qquad \leq \eltwo{e^{-A\prn{t-t_k}}\prn{a_{\Gk}(t_k) - \adag}} \\
 & \qquad \qquad \qquad \qquad + \eltwo{\prn{I-e^{-A\prn{t-t_k}}}\prn{a_{\Gk}^{\infty}-\adag}} \\
 & \qquad \leq e^{-\prn{1-\delta}\prn{t-t_k}} \eltwo{a_{\Gk}(t_k)-\adag} \\
 & \qquad \qquad \qquad \qquad + \prn{1-e^{-\prn{1+\delta}\prn{t-t_k}}} \eltwo{a_{\Gk}^{\infty}-\adag} \\
 & \qquad \leq e^{-\prn{1-\delta}\prn{t-t_k}} C_{\delta}(p) + \prn{1-e^{-\prn{1+\delta}\prn{t-t_k}}} \dfrac{1-\delta}{1+\delta}C_{\delta}(p) \\
 & \qquad \overset{(i)}{\leq} C_{\delta}(p). 
\end{align*}
To prove the last inequality, we study the function:
$$h(t) = \prn{1-e^{-(1-\delta)t}} C - \prn{1-e^{-\prn{1+\delta}t}} \dfrac{1-\delta}{1+\delta} C.$$
The derivative of $h(t)$ is
\begin{align*}
 h'(t) & = \prn{1-\delta}e^{-\prn{1-\delta}t}C - \prn{1-\delta}e^{-\prn{1+\delta}t}C \\
 & = \prn{1-\delta}C\bigprn{e^{-\prn{1-\delta}t} - e^{-\prn{1+\delta}t}} \geq 0.
\end{align*}
Since $h(0) = 0$, and $h'(t)\geq0$ for all $t\geq0$, we have $h(t)\geq0$ for all $t\geq0$, and the inequality (i) holds. 

Finally, since the vector $a(t)-\adag$ is continuous with time:
\begin{equation*}
\eltwo{a_{\G_{k+1}}(t_{k+1}) - \adag} = \eltwo{a_{\Gk}(t_{k+1}) - \adag}  \leq C_{\delta}(p). \tag*{\IEEEQED}
\end{equation*}
\let\IEEEQED\relax 
\end{IEEEproof}
%

\section{Proof of Theorem~\ref{th:opt}}
\label{app:opt}
\begin{IEEEproof}
To prove that the active set $\G$ is a subset of $\G_{\dagger}$ for all time $t\geq0$, we show by induction that for all switching times $t_k$, and for all time $t\in\prn{t_k,t_{k+1}}$, we have:
\begin{equation}
\abs{u_j(t)} \leq \lambda, \qquad \qquad \forall j\in\Gc_{\dagger}.
\label{eq:inactive3}
\end{equation}
If this condition is satisfied, then nodes in $\Gc_{\dagger}$ stay below threshold and the next active set $\G_{k+1}$ is also a subset of $\G_{\dagger}$. We add a second induction hypothesis:
\begin{equation}
\eltwo{a_{\G_k}(t)-\adag} \leq C_{\delta}(S) \qquad \forall t\in\Prn{t_k,t_{k+1}}.
\label{eq:rcond}
\end{equation}
By the theorem hypotheses, the initial active set is a subset of $\G_{\dagger}$ and \eqref{eq:icond} holds, so \eqref{eq:inactive3} and \eqref{eq:rcond} hold at $t=0$. We now assume that the two induction hypotheses hold for a particular switching time $t_k$. 
If there is no more switching after $t_k$, then we are done. Otherwise, using the dynamics in \eqref{eq:inactsol1}, we know that, for all $j\in\Gc_{\dagger}\subset\Gc_k$, we have $\forall t\in\Prn{t_k,t_{k+1}}$:
$$u_j(t) = e^{-(t-t_k)} u_j^{t_k} + e^{-t} \dint{t_k}{t}{e^{\nu} \rho_j(\nu) d\nu},$$
with $\rho_j(\nu) = \Phi_j^T\prn{y-\Phi_{\G_k}a_{\G_k}(\nu)}.$
We bound the absolute value of the expression above using:
\begin{align*}
\abs{u_j(t)} & = \abs{e^{-(t-t_k)} u_j^{t_k} + e^{-t} \dint{t_k}{t}{e^{\nu} \rho_j(\nu) d\nu}} \\
 & \leq e^{-(t-t_k)} \abs{u_j^{t_k}} + e^{-t} \dint{t_k}{t}{e^{\nu}\abs{\rho_j(\nu)}d\nu} \\
 & \leq e^{-(t-t_k)} \abs{u_j^{t_k}} + \prn{1-e^{-(t-t_k)}} \sup_{\nu'\in\Prn{t_k,t_{k+1}}}\abs{\rho_j(\nu')}.
\end{align*} 
Since at time $t_k$, node $j\in\Gc_{\dagger}$ is inactive, we have: $\abs{u_j^{t_k}} \leq\lambda$. As a consequence, condition \eqref{eq:inactive3} is satisfied if:
\begin{equation}
\sup_{\nu'\in\Prn{t_k,t_{k+1}}}\abs{\rho_j(\nu')} \leq \lambda.
\label{eq:bound}
\end{equation}
%
We use the fact that the matrix $\Phi_j^T\Phi_{\G_{\dagger}}$ is a submatrix of $\Phi^T\Phi-I$ with $(S+1)$ distinct columns and apply the RIP of order $(S+1)$: $\norm{\Phi_j^T\Phi_{\G_{\dagger}}}\leq \delta$. Then, we have that for all time $t\in\Prn{t_k,t_{k+1}}$ and for all nodes $j\in\Gc_{\dagger}$:
\begin{align*}
\abs{\rho_j(t)} & = \abs{\Phi_j^T\prn{y-\Phi_{\G_k}a_{\G_k}(t)}} \\
 & = \abs{\Phi_j^T\prn{\Phi_{\G_{\dagger}}\adag + \epsilon-\Phi_{\G_k}a_{\G_k}(t)}} \tag{$y=\Phi_{\G_{\dagger}}\adag+\epsilon$}\\
 & = \abs{\Phi_j^T\Phi_{\G_{\dagger}} \prn{\adag - a_{\G_k}(t)} +  \Phi_j^T\epsilon}  \tag{since $\G_k\subset\G_{\dagger}$}\\
 & \leq \abs{\Phi_j^T\Phi_{\G_{\dagger}} \prn{\adag - a_{\G_k}(t)}} + \abs{\Phi_j^T\epsilon} \\
 & \leq \norm{\Phi_j^T\Phi_{\G_{\dagger}}}\eltwo{\adag - a_{\G_k}(t)} + \elinf{\Phi_{\Gc_{\dagger}}^T\epsilon} \\
 & \leq \delta\eltwo{\adag - a_{\G_k}(t)} +  \elinf{\Phi_{\Gc_{\dagger}}^T\epsilon}.
\end{align*}
%
%
We apply Lemma~\ref{lem:timebound} to get a bound that holds uniformly across time: $\eltwo{\adag - a(t)}\leq C_{\delta}(S)$,  $\forall t\in\Prn{t_k,t_k+1}$. In particular $ \eltwo{a_{\G_{k+1}}(t_{k+1}) - \adag} \leq C_{\delta}(S)$ and the induction hypothesis \eqref{eq:rcond} remains true at time $t_{k+1}$. 

Putting the pieces together and using condition \eqref{eq:rel1}, we have for all time $t\in\Prn{t_k,t_{k+1}}$ and for all nodes $j\in\Gc_{\dagger}$:
\begin{align*}
 \abs{\rho_j(t)} & \leq \delta C_{\delta}(S) + \elinf{\Phi_{\Gc_{\dagger}}^T\epsilon}  \leq \lambda \prn{1-\alpha\delta\sqrt{S} + \alpha\delta\sqrt{S}} = \lambda
\end{align*}
This shows that \eqref{eq:inactive3} holds for all time $t\in\Prn{t_k,t_{k+1}}$, which ends the proof by induction of the theorem.
\end{IEEEproof}


\section{Proof of Theorem~\ref{th:size}}
\label{app:size}

\begin{IEEEproof}
We want to show that for all time during convergence, no more than $q$ nodes are active, i.e. $\abs{\G} \leq q$. First, we introduce some notations and denote by $\GD(t)$ the set containing the $q$ biggest nodes in $u(t)$. This set depends on time, but we will often remove the dependence in the notation for readability. By definition of $\GD(t)$, we have $\forall j\in\GD^c(t)$
\begin{equation}
\abs{u_j(t)} \leq { \eltwo{u_{\GD(t)}(t)}}/{\sqrt{q}}.
\label{eq:DGCbound}
\end{equation}
The theorem holds if for all time $t\geq0$, nodes $j$ in $\GD^c(t)$ are below threshold: $\abs{u_j(t)} \leq \lambda$. Thus, we prove the theorem by induction, by showing that for all switching times $t_k$, we have $\forall t\in\Prn{t_k,t_{k+1}}$
\begin{equation}
\eltwo{u_{\GD(t)}(t)} \leq \lambda\sqrt{q},
\label{eq:rcond2}
\end{equation}
along with the following induction hypothesis:
\begin{equation}
 \eltwo{ a(t) - \adag} \leq C_{\bar{\delta}}(q)
 \label{eq:rcond3}
\end{equation}
By \eqref{eq:initcond2}, \eqref{eq:rcond2} holds at $t=0$ and
\begin{align*}
 \eltwo{a(0) - \adag} & \leq  \eltwo{a(0)} + \eltwo{\adag}\leq \eltwo{u(0)} + \eltwo{\adag} \\
 & \leq \lambda\sqrt{q} + \eltwo{\adag} \; \leq C_{\bar{\delta}}(q),
\end{align*}
so \eqref{eq:rcond3} also holds at $t=0$.

We now assume that for some switching time $t_k$, \eqref{eq:rcond2} and \eqref{eq:rcond3} hold. If there is no more switching, we are done. Otherwise, by \eqref{eq:odesol2}, the dynamics on $\GD$ for $t\in\prn{t_k,t_{k+1}}$ are: 
$$ u_{\GD}(t) = e^{-(t-t_k)}u_{\GD}(t_k) + e^{-t}\dint{t_k}{t}{ e^{\nu} \rho_{\GD}(\nu) d\nu},$$
where $\rho_{\GD}(\nu) = a_{\GD}(\nu) - \Phi_{\GD}^T\Phi a(\nu) + \Phi_{\GD}^Ty$. We bound the $\ell_2$-norm of this quantity as follows:
\begin{align}
 & \eltwo{u_{\GD}(t)} \leq e^{-(t-t_k)}\eltwo{u_{\GD}(t_k)} \nonumber \\
 & \qquad \qquad \qquad + e^{-t}\dint{t_k}{t}{ e^{\nu} \sup_{\nu'\in{t_k,t_{k+1}}}\eltwo{\rho_{\GD}(\nu')} d\nu} \nonumber \\
 & ~ \leq e^{-(t-t_k)} \eltwo{u_{\GD}(t_k)} + \prn{1-e^{-(t-t_k)}} \sup_{\nu'\in{t_k,t_{k+1}}}\eltwo{\rho_{\GD}(\nu')}. \label{eq:uGDbound}
\end{align}
By the induction hypothesis \eqref{eq:rcond2}, we know that
$$\eltwo{u_{\GD}(t_k)} \leq \lambda\sqrt{q}.$$
We now find a bound for all $t\in\prn{t_k,t_{k+1}}$ for:
\begin{align*}
 & \eltwo{\rho_{\GD}(t)} = \eltwo{ a_{\GD}(t) - \Phi_{\GD}^T\Phi a(t) + \Phi_{\GD}^Ty} \\
 & ~~ = \eltwo{\adag + \prn{I - \Phi_{\GD}^T\Phi}(a(t) - \adag) + \Phi_{\GD}^T\epsilon} \\
 & ~~ \leq \eltwo{\adag} + \norm{I - \Phi_{\GD}^T\Phi_{\G_{\dagger}\cup\Gk}}\eltwo{a(t) - \adag} + \eltwo{\Phi_{\GD}^T\epsilon}.
\end{align*}
By \eqref{eq:rcond2} and \eqref{eq:DGCbound}, $\GD^c \subset\Gk^c$, so $\Gk\subset\GD$. Thus, we can apply the RIP of order $S+p$ to the matrices \mbox{$I - \Phi_{\GD}^T\Phi_{\G_{\dagger}\cup\GD}$} and $\Phi_{\GD}$. Moreover, since \eqref{eq:rcond3} holds at time $t_k$, we apply Lemma~\ref{lem:timebound}, which proves that the induction condition \eqref{eq:rcond3} is true at time $t_{k+1}$ and gives a uniform bound on the quantity $\eltwo{a(t) - \adag}$:
\begin{align*}
 & \eltwo{\rho_{\GD}(t)} \leq \eltwo{\adag} + \bar{\delta}C_{\bar{\delta}}(q) + (1+\bar{\delta})\eltwo{\epsilon} \\
 & ~~ = \prn{1+\bar{\delta}(1+\bar{\delta})(1-\bar{\delta})^{-2}}\eltwo{\adag} \\
 & \qquad \qquad + (1+\bar{\delta})\prn{\bar{\delta}(1-\bar{\delta})^{-2}\sqrt{1-\bar{\delta}} +1}\eltwo{\epsilon} \\
 & \qquad \qquad + \bar{\delta}(1+\bar{\delta})(1-\bar{\delta})^{-2}\lambda\sqrt{q} \\
 & ~~ \leq (1+\bar{\delta})(1-\bar{\delta})^{-2} \prn {\eltwo{\adag} + \sqrt{1-\bar{\delta}}\eltwo{\epsilon} + \bar{\delta}\lambda\sqrt{q}}.
\end{align*}
Applying the theorem hypothesis \eqref{eq:rel2}, we get
\begin{align*}
 \eltwo{\rho_{\GD}(t)} & < (1-\bar{\delta})^{-2}\prn{1-3\bar{\delta} + \bar{\delta}(1+\bar{\delta})}\lambda\sqrt{q} = \lambda\sqrt{q}.
\end{align*}
Plugging this back into \eqref{eq:uGDbound} shows that $\eltwo{u_{\GD}(t)} \leq \lambda\sqrt{q}$ for all $t\in\Prn{t_k,t_{k+1}}$. In particular, we proved that the induction condition \eqref{eq:rcond2} also holds, which finishes the proof.
\end{IEEEproof}




\bibliographystyle{IEEEbib}
\bibliography{LcaCS}

\end{document}